\documentclass[12pt, reqno]{amsart}
\usepackage{amsmath, amsthm, amscd, amsfonts, amssymb, graphicx, color, mathrsfs}
\usepackage[bookmarksnumbered, colorlinks, plainpages]{hyperref}
\usepackage[all]{xy}
\usepackage{slashed}

\newtheorem{theorem}{Theorem}[section]

\newtheorem{corollary}[theorem]{Corollary}
\theoremstyle{definition}

\theoremstyle{remark}

\numberwithin{equation}{section}

\allowdisplaybreaks

\begin{document}
	\setcounter{page}{1}
	
\title[Spectral Multipliers of the Anharmonic Oscillator]{$L^p$-$L^q$ Boundedness of Spectral Multipliers of the Anharmonic Oscillator\\
$L^p$-$L^q$ continuit\'e des multiplicateurs spectraux de l'oscillateur anharmonique}

	\author[Marianna Chatzakou]{Marianna Chatzakou}
	\address{
		Marianna Chatzakou:
		\endgraf
		Department of Mathematics: Analysis Logic and Discrete Mathematics
		\endgraf
		Ghent University, Belgium
		\endgraf
		{\it E-mail address} {\rm Marianna.Chatzakou@UGent.be}
		\endgraf
	}
	
	\author[Vishvesh Kumar]{Vishvesh Kumar}
	\address{
		Vishvesh Kumar:
		\endgraf
		Department of Mathematics: Analysis Logic and Discrete Mathematics
		\endgraf
		Ghent University, Belgium
		\endgraf
		{\it E-mail address} {\rm vishveshmishra@gmail.com}
		\endgraf
	}



	
	
	
	\begin{abstract} In this note we study the $L^p-L^q$ boundedness of Fourier multipliers of anharmonic oscillators, and as a consequence also of spectral multipliers, for the range $1<p \leq 2 \leq q <\infty$. The underlying Fourier analysis is associated with the eigenfunctions of an anharmonic oscillator in some family of differential operators having derivatives of any order. Our analysis relies on a version of the classical Paley-type inequality, introduced by H\"ormander, that we extend in our nonharmonic setting.\\
	{\bf R\'esum\'e}\\
	Dans cette note, nous étudions la $L^p$-$L^q$ continuité des multiplicateurs de Fourier des oscillateurs anharmoniques, et par conséquent des multiplicateurs spectraux également, pour $1<p \leq 2 \leq q <\infty$. L'analyse de Fourier sous-jacente est associée aux fonctions propres d'un oscillateur anharmonique dans certaines familles d'opérateurs différentiels ayant des dérivées d'ordre quelconque. Notre analyse s'appuie sur une version de l'inégalité classique de type Paley, introduite par Hörmander, que nous étendons dans notre cadre non harmonique.
\end{abstract} \maketitle
	
	
\section{Introduction}
The aim of this note is to study the $L^p-L^q$ boundedness conditions for the Fourier and spectral multipliers where the pertinent Fourier analysis is developed  with respect to a member of the family of anharmonic oscillators on $\mathbb{R}^n$. Such operators were considered in \cite{CDR18} as ``prototype'' for the spectral analysis that is performed there. 

\par The boundedness of Fourier multipliers is a central topic in harmonic analysis and its origin can be traced back to H\"ormander's seminal paper \cite{Hormander1960} in 1960. The results concerning to the boundedness of multipliers are of $L^p$-type, or $L^p-L^q$-type. In the first case, Mihlin-H$\ddot{\text{o}}$rmander \cite{M} or Marcinkiewicz theorems \cite{CS} come into play providing conditions on the regularity of the symbol, while in the second case, the sharp decay of the spectral information is analysed. 
\par Regarding the $L^p-L^q$ boundedness that we study here, let us recall H\"ormander's Fourier multiplier theorem as in \cite{Hormander1960}: For $1<p\leq 2 \leq q <\infty,$  the Fourier multiplier $T_\sigma: \mathcal{S}(\mathbb{R}^n) \rightarrow \mathcal{S}'(\mathbb{R}^n)$ associated with symbol $\sigma:\mathbb{R}^n \rightarrow \mathbb{C}$  defined by \[\mathcal{F}(T_\sigma f)(\xi)= \sigma(\xi) \mathcal{F}(f)(\xi)\] for $\xi \in \mathbb{R}^n$, where $\mathcal{F}$ denotes the Euclidean Fourier transform of $f$, has a bounded extension from $L^p(\mathbb{R}^n)$ to $L^q(\mathbb{R}^n)$ provided that the symbol satisfies the condition
\begin{equation}\label{Horm.thm}
   \sup_{s>0} s \left( \int_{\xi \in \mathbb{R}^n: |\sigma(\xi)| \geq s}d\xi \right)^{\frac{1}{p}-\frac{1}{q}}<+\infty\,.
\end{equation} 
\par The $L^p-L^q$ boundedness of spectral and Fourier mutlipliers has been studied in different settings; we refer for instance to the works \cite{CGM} on symmetric spaces, and \cite{AR,ARN,ARN2} on groups. In \cite{AR,ARN2} the authors suggest an approach that is also effective for our purposes: to establish the Hausdorff-Young-Paley inequality by interpolating between the Hausdorff-Young inequality and Paley-type inequality. However, the machinery we use here differs from the one in \cite{AR} and relies on the underlying Fourier analysis, and so on the spectral properties of our ``model operator'', and on the Sobolev estimates associated with it as acquired in \cite{CDR18}.  
\par Precisely, the Fourier analysis we present here is the natural analogue of the one developed in \cite{Ruz-Tok} in the general setting of nonharmonic analysis associated with a ``model operator''. A concrete example of the latter is our choice of the anharmonic oscillator with a potential that is a homogeneous polynomial of a certain degree, and with derivatives of any order; in particular the operator we consider is of the form \footnote{The addition of the constant $1$ to the potential is required for the invertibility of the operator and for controlling the lower bound of the spectrum of the operator.}
\[
(-\Delta)^{l}+|x|^{2k}+1\,, \quad x \in \mathbb{R}^n\,,
\]
where $l,k \geq1$ are some fixed integers, $|\cdot|$ stands for the Euclidean norm on $\mathbb{R}^n$, and $\Delta$ is the usual Laplace operator on $\mathbb{R}^n$. Under such analysis a Fourier multiplier is $L^p-L^q$ bounded if a condition similar to \eqref{Horm.thm} is satisfied for a symbol of the multiplier adapted to our setting.

\section{Notation and preliminaries}
In order to present our results, let us briefly recall the necessary notation and concepts involved in our analysis. For a more detailed discussion on concepts related to the spectral analysis of the anharmonic oscillator we refer to \cite{CDR18,chatzakou2020lplq}, and for those related to the associated Fourier analysis we refer to \cite{chatzakou2020lplq}, or, for a more general consideration, to \cite{Ruz-Tok}.
\subsection{The anharmonic oscillator and its spectrum}
The \textit{anharmonic oscillator} is a differential operator on $\mathbb{R}^n$ that is served as a generalisation of the harmonic oscillator. Here we consider an anharmonic oscillator in the family of operators $(A_{k,l})_{k,l}$, for different values of the integer parameters $k,l\geq 1$, defined as follows
\begin{equation*}\label{defn.anh}
    A_{k,l}:=(-\Delta)^{l}+|x|^{2k}+1\,,\quad x \in \mathbb{R}^n\,.
\end{equation*}
Here and thereafter we assume that the parameters $k,l$ are fixed, and the notation will be simplified to $A:=A_{k,l}$.

\par The operator $A$ is densely defined on its domain $\mathcal{S}(\mathbb{R}^n)\subset L^2(\mathbb{R}^n)$ and admits a (unique) self-adjoint extension on $L^2(\mathbb{R}^n)$. Its spectrum is purely discrete and lies in the interval $[1,\infty)$. If we keep the same notation for its self adjoint extension, then by the functional calculus we have, for $m \in \mathbb{R},$
\begin{equation*}\label{dom.Ar}
Dom\left(A^{m}\right)= \left\{  u \in L^2(\mathbb{R}^n) : \sum_{j=1}^{\infty} {\lambda_j}^{2m} |\langle u_j,u \rangle_{L^2}|^2 <\infty \right\},
\end{equation*}
where the set $(u_j)_j$, with $\|u_j\|_{L^2}=1$, is the set of normalized eigenfunctions of $A$, and $(\lambda_j)_j$ is the corresponding set of eigenvalues. Thus, for $u \in Dom \left(A^{m}\right)$, we have
\[
A^{m}u=\sum_{j=1}^{\infty} \lambda_{j}^{m} \langle u_{j},u \rangle_{L^2} u_j\,.
\] 
Despite the intensive research on the eigenvalue problem associated with the anhramonic oscillator, the exact solution of it, even in the quartic case, is still explicitly unknown. We refer to \cite{CDR18} for a recent exposition of the related works. Let us summarise the spectral information that is known in the case of the anharmonic oscillator $A$: For the inverse operator $A^{-1}$ we know that 
  $A^{-1} \in S_r(L^2(\mathbb{R}^n))$ for $r>\frac{(k+l)n}{2kl}$, where by $S_r(L^2(\mathbb{R}^n))$ we have denoted the \textit{$r$th Schatten-von Neumann class} of operators on the Hilbert space $L^2(\mathbb{R}^n)$, see  \cite[Corollary 5.3]{CDR18}.
 As a consequence, for the eigenvalue counting function $N(\lambda):= \# \{  j: \lambda_j < \lambda\},$ we get the following asymptotic estimate 
 \begin{equation}
    \label{eig.count.f}
    N(\lambda) \lesssim \lambda^{n\left(\frac{1}{2k}+\frac{1}{2l}\right)}\,,\quad \text{as}\quad \lambda \rightarrow \infty\,,
\end{equation}
where we have assumed that the eigenvalues $(\lambda_j)_j$ are arranged in the decreasing order; see \cite[Remark 5.7]{CDR18}, or \cite[Theorem 5.5]{HR}.

\subsection{Fourier analysis associated with the anharmonic oscillator}
For the setup of this section we follow \cite{Ruz-Tok}.\\
\par The space of functions 
\begin{equation*}\label{Eq:R-T-dom}
C^\infty_{A}(\mathbb{R}^n):=\cap_{m=1}^\infty \textnormal{Dom}(A^m)\,,
\end{equation*}
will be the domain of the Fourier operator $\mathcal{F}_A$ in our setting. The space $C^{\infty}_{A}$ is a Fr\'echet space if endowed with the topology induced by the family of norms  
$$\|f\|_{C^m_{A}}:=\max_{k\leq m} \| A^k f\|_{L^2(\mathbb{R}^n)},\,\, m\in\mathbb{N}_0,\,\,f\in  C^\infty_{A}(\mathbb{R}^n).$$
Then, defining the \textit{$A$-Fourier transform} by \[
(\mathcal{F}_A f)(j):=\int_{\mathbb{R}^n}f(x)\overline{u_j(x)}\,dx\,,
\]
we see that $\mathcal{F}_A(C^{\infty}_{A})\subset \mathcal{S}(\mathbb{N})$, where $\mathcal{S}(\mathbb{N})$ is the space of rapidly decreasing functions $\phi:\mathbb{N}\to \mathbb{C}$, i.e., for any $N<\infty,$ there exists a constant $C_{\phi,N}>0$ such that \[|\phi(j)|\leq C_{\phi,N}\langle j\rangle^{-N}\,\,\text{for all } j\in\mathbb{N}\,.\] The Fr\'echet topology of $\mathcal{S}(\mathbb{N})$ is induced by the family of semi-norms $\{p_r\}_{r \in \mathbb{N}}$ given by \[p_r(\phi):=\sup_{j\in\mathbb{N}} \langle j\rangle^r |\phi(j)|\,.\]
With the topological structures as described above, the map $\mathcal{F}_A$ becomes a bijective homeomorphism form $C^{\infty}_{A}$ to $\mathcal{S}(\mathbb{N})$, and the Fourier inversion formula becomes available in our setting. Moreover, for $f \in L^2$, the operator $\mathcal{F}_A$ becomes a bijective isomorphism onto the Hilbert space $l^{2}_{A}$ defined as the space of sequences $f=(f_j)_{j \in \mathbb{N}}$ under the norm
\[
\|f\|_{l^{2}_{A}}:=\left(\sum_{j \in \mathbb{N}}|f_j|^2 \right)^{\frac{1}{2}}\,.
\]
The latter reads as the \textit{Plancherel identity} in our setting, i.e., in symbols: $\|f\|_{L^2}=\|\mathcal{F}_A(f)\|_{l^{2}_{A}}$.






\section{Main results}
The Hausdorff-Young inequality in our setting is derived verbatim as a particular case of the Fourier analysis associated with a ``model operator'' \cite{Ruz-Tok}. The Paley-type inequality is more involved and its proof can be found in the preprint \cite[Theorem 4.1]{chatzakou2020lplq}. Then, by means of the Stein-Weiss interpolation theorem, we can extend the Hausdorff-Young-Paley inequality in our setting \cite[Theorem 4.4]{chatzakou2020lplq} which is the key for proving our result on the boundedness of Fourier multipliers as follows. We refer to \cite{chatzakou2020lplq} for the detailed proofs of the results below.

\par First, let us recall that $T_\sigma$ is an $A$-Fourier multiplier, if it satisfies the identity
$$\mathcal{F}_{A}({T_\sigma})(j)= \sigma(j) \mathcal{F}_{A}{f}(j),\,\,\,\, j \in \mathbb{N},$$
 where $\sigma :\mathbb{N} \rightarrow \mathbb{C}$ is a function that can be served as the symbol of the multiplier $T_\sigma$ for our analysis.

\begin{theorem}\label{Th:LpLq-1}
Let $1<p \leq 2 \leq q <\infty.$  Let $T_\sigma$ be a $A$-Fourier multiplier with symbol $\sigma$. Then we have 
\begin{equation}\label{thm1.eq} \|T_\sigma\|_{L^p(\mathbb{R}^n) \rightarrow L^q(\mathbb{R}^n)} \lesssim \sup_{s>0} s\left(  \ \sum_{\overset{j \in \mathbb{N}}{|\sigma(j)| \geq s}} \|u_j\|^2_{L^\infty(\mathbb{R}^n)}  \right)^{\frac{1}{p}-\frac{1}{q}}\,.
\end{equation}
\end{theorem}
Let us mention that \eqref{thm1.eq} is well-defined since for the $L^\infty$-estimate of eigenfunctions we have: $\|u_j\|_{L^\infty}\leq C \sqrt{\lambda_j}$, for some $C>0$ and for all $j$; see \cite[Lemma 2.2]{chatzakou2020lplq}. We note here that Theorem \ref{Th:LpLq-1} can be thought as a particular case of recently obtained result in \cite{DKRN20} for a more general setting.

\par As a consequence of Theorem \ref{Th:LpLq-1} and of the spectral decay of the anharmonic oscillator $A$ in \eqref{eig.count.f} we derive the following result on the boundedness of the spectral multipliers in our setting.

\begin{theorem}\label{cor.anh.1}
  Let $1<p\leq 2 \leq q<\infty$, and let $\varphi$ be any decreasing function on $[1,\infty)$ such that
 $\lim_{u \rightarrow \infty} \varphi(u)=0$. Then,
  \begin{equation*}
  \label{cor.1.estim}
  \|\varphi(A)\|_{L^p(\mathbb{R}^n) \rightarrow L^q(\mathbb{R}^n)} \lesssim \sup_{u\geq 1} \varphi(u)\left(u^{1+\frac{(k+l)n}{2kl}}\right)^{\frac{1}{p}-\frac{1}{q}}\,.
  \end{equation*}
 
 \end{theorem}

 \par  Let us summarise some applications of Theorem \ref{cor.anh.1}.
 \begin{corollary}\label{cor1}
 Let $1<p\leq 2 \leq q < \infty$. Then, the operator $A^{-m}$ is bounded from $L^p(\mathbb{R}^n)$ to $L^q(\mathbb{R}^n)$  for $m\geq  \left(1+\frac{(k+l)n}{2kl}\right) \left(\frac{1}{p}-\frac{1}{q}\right)$. 
 \end{corollary}
 Corollary \ref{cor1} implies in particular that $A^{-m}$ is bounded from $L^2(\mathbb{R}^n)$ to $L^2(\mathbb{R}^n)$ for any $m \geq 0$. This is known in more generality by the general theory of pseudo-differential operators, see Beals \cite[Theorem 3.1]{B}, which gives us $A^{-m}\in \mathcal{L}(L^p)$, for any $m>0$, and for any $p \in (1,\infty)$. However, using Corollary \ref{cor1} we can study the $L^p$-$L^q$ boundedness for $p \neq q$.
 
 \begin{corollary}\label{cor2}
 If $v=v(x,t)$ is the solution to the heat equation associated to $A$, i.e., if we have
     	\begin{equation*}
	\label{heat.eq}
	\partial_{t}v+Av=0\,,\quad v(0)=v_0\,,
	\end{equation*}
	then the time decay rate for $v$ is estimated as:
     \[
	\|v(t,\cdot)\|_{L^q(\mathbb{R}^n)} \lesssim t^{-\left( 1+\frac{(k+l)n}{2kl}\right)\left(\frac{1}{p}-\frac{1}{q} \right)}\ \|v_0\|_{L^q(\mathbb{R}^n)}\,,
	\]
	where $1<p\leq 2 \leq q < \infty$.
 \end{corollary}
 \begin{corollary}\label{cor3}
 Let $1<p\leq 2 \leq q < \infty$. Then, for $m\geq \left(1+\frac{(k+l)n}{2kl} \right)\left(\frac{1}{p}-\frac{1}{q} \right)$ we obtain the Sobolev-type estimate:
 	\[
		\|f\|_{L^q(\mathbb{R}^n)}\lesssim  \|A^{m}f\|_{L^p(\mathbb{R}^n)}\,.
		\]	
 \end{corollary}
 
\section*{acknowledgement}
 The authors are supported  by the FWO  Odysseus  1  grant  G.0H94.18N:  Analysis  and  Partial Differential Equations and by the Methusalem programme of the Ghent University Special Research Fund (BOF)
(Grant number 01M01021).

We thank Michael Ruzhansky and Duvan Cardona for helpful discussion and suggestions.

\bibliographystyle{amsplain}

\end{document}